\documentclass[reqno]{amsart}
\usepackage{url}
\newtheorem{theorem}{Theorem}[section]

\theoremstyle{definition}
\newtheorem{definition}[theorem]{Definition}

\theoremstyle{remark}

\numberwithin{equation}{section}
\usepackage{algorithmic}
\usepackage{amssymb}
\usepackage[top=3cm,bottom=2cm,right=2cm,left=2cm]{geometry}
\usepackage{hyperref}

\begin{document}

\title[A finite-bound partition equinumerosity result]{A finite-bound
  partition equinumerosity result generalizing\\
  a solution of a problem posed by Andrews and Deutsch}

\author{Michael J.\ Schlosser}
\address{Fakult\"{a}t f\"{u}r Mathematik, Universit\"{a}t Wien\\
Oskar-Morgenstern-Platz 1, A-1090, Vienna, Austria}
\email{michael.schlosser@univie.ac.at}
\thanks{The first author was partially
  supported by Austrian Science Fund FWF
  \href{https://doi.org/10.55776/P32305}{10.55776/P32305}.}

\author{Nicolas Allen Smoot}
\address{Fakult\"{a}t f\"{u}r Mathematik, Universit\"{a}t Wien\\
Oskar-Morgenstern-Platz 1, A-1090, Vienna, Austria}
\email{nicolas.allen.smoot@univie.ac.at}
\thanks{The second author was fully
  supported by Austrian Science Fund FWF
  \href{https://doi.org/10.55776/PAT6428623}{10.55776/PAT6428623}.}

\keywords{integer partitions; partition bijection}

\subjclass[2020]{Primary 05A19, Secondary 11P81}

\date{}

\begin{abstract}
We introduce a finite-bound extension of a partition equinumerosity result which was orignally proposed as a problem by Andrews and Deutsch in 2016, and given a generalized form in 2018 by Smoot and Yang.  We also give a simple bijective proof, itself an extension of the proof given of the 2018 generalization.
\end{abstract}

\maketitle

\section{Background}

In 2016 Andrews and Deutsch proposed the following problem:

\noindent [American Mathematical Monthly Problem 11908]\ \textit{Let $n, k$ be nonnegative integers.  Show that the number of partitions of $n$ having $k$ even parts is the same as the number of partitions of $n$ in which the largest repeated part is $k$ (defined to be 0 if the parts are all distinct).}

As an example, they list the three partitions of 7 with two even parts:
\begin{align*}
&4+2+1,\ \ \ \ \ 3+2+2,\ \ \ \ \ 2+2+1+1+1,
\end{align*} and the three partitions of 7 in which the largest repeated part is 2:
\begin{align*}
&3+2+2,\ \ \ \ \ 2+2+2+1,\ \ \ \ \ 2+2+1+1+1.
\end{align*}  The counting functions for each class of partitions has a natural generating function, and Madhyastha gave a proof of the problem through manipulation of their generating functions \cite[Solution I]{SmootA}, showing that
\begin{align}
\frac{q^{2k}}{(q^2;q^2)_k} \frac{(q^2;q^2)_{\infty}}{(q;q)_{\infty}} = \frac{q^{2k}}{(q;q)_k}\frac{(q^{2(k+1)};q^2)_{\infty}}{(q^{k+1};q)_{\infty}}.\label{solutionI}
\end{align}
Here, and in the first proof of Theorem~\ref{stheorem}, $q$ is a formal
variable; for any indeterminate $a$ and non-negative integer $k$
(we also allow $k=\infty$) the $q$-shifted factorial is defined as
$(a;q)_k=\prod_{0\le j< k}(1-aq^j)$ (empty products being $1$).

However, two additional discoveries have been made regarding this problem.  First, it may be generalized, as was discovered by Nicolas Allen Smoot and Mingjia Yang \cite[Solution II, Editorial comment]{SmootA}, \cite[Theorem 23]{Johnson}:
\begin{theorem}[N. Smoot, M. Yang]\label{sythm}
Let $n, k, d$ be nonnegative integers.  The number of partitions of $n$ with $k$ parts divisible by $d$ is the same as the number of partitions of $n$ in which $k$ is the largest part to occur at least $d$ times.
\end{theorem}
Moreover, a delightful proof by direct construction may be given for this generalization \cite[Solution II]{SmootA}.
\begin{proof}
Let $\lambda$ be a partition of $n$ with $k$ parts divisible by $d$.  Write
\begin{align*}
\lambda = \mu + o,
\end{align*} in which $\mu$ is the subpartition containing the $k$ parts divisible by $d$, and $o$ contains the remaining parts.

We map $\mu$ to its conjugate, which we denote $\mu^{\ast}$.  Notice that $\mu^{\ast}$ contains $k$ as its largest part, and that indeed $k$ must repeat $dq_k$ times for some integer $q_k > 0$.

Next we map $o$, which contains no parts divisible by $d$, to its corresponding partition $\delta$ in which every part may only repeat up to $d-1$ times---say, using Glaisher's map \cite{Glaisher}.  Both of these actions constitute bijections.  Our new partition, then, is 
\begin{align*}
\kappa = \mu^{\ast} + \delta.
\end{align*}  Notice that any part $i$ of $\kappa$ which is less than or equal to $k$ will occur $dq_i+r_i$ times, in which $dq_i$ comes from $\mu^{\ast}$ and $q_i$ is a nonnegative integer (which is possibly 0 except when $i=k$), and $r_i$, with $0\le r_i\le d-1$ comes from $\delta$.  Thus parts less than $k$ may repeat any number of times.  On the other hand, any part greater than $k$ must only have a contribution from $\delta$, and must therefore occur strictly less than $d$ times.

This mapping is well-defined, as well as easily reversable: for any partition $\kappa$ in which $k$ is the largest part to occur at least $d$ times, take all parts $i$ which occur $dq_i+r_i$ times with $q_i \ge 0$ (which is positive for $i=k$), and define the subpartition $\mu^{\ast}$ by collecting $i$ to occur $dq_i$ times.  The largest part will be $k$.  What is left will be parts $i$ which repeat $r_i$ times for $0\le r_i < d$.  Collect these terms into $\delta$.

We map $\mu^{\ast}$ to its conjugate $\mu$, which will contain $k$ parts, all divisible by $d$.  Also we map $\delta$ to $o$, which contains no multiples of $d$.  Our new partition is $\lambda=\mu + o$.
\end{proof}

\section{A finite-bound generalization of Theorem~\ref{sythm}}

The proof of Theorem~\ref{sythm} is especially satisfying, but it is not obvious that anything more can be done with it.  However, we have found that Theorem \ref{sythm} may be extended to a finite-bound result.
\begin{definition}
Let $A(n,k,d,m)$ be the set of integer partitions of $n$ in which that exactly $k$ parts are divisible by $d$, and all other parts are strictly less than $md$.
\end{definition}
\begin{definition}
Let $B(n,k,d,m)$ be the set of integer partitions of $n$ such that:
\begin{itemize}
\item For $m < k$, the largest part is $kd$, and all parts greater than $md$ are divisible by $d$;
\item For $m\ge k$, the part $k$ occurs at least $d$ times, none of the parts exceed $md$, and any part $i$ such that $k < i\le m$ occurs less than $d$ times.
\end{itemize}
\end{definition}\noindent With these definitions, we have the following:
\begin{theorem}\label{stheorem}
For $n,k,d,m\in\mathbb{Z}_{\ge 1}$, $\left|A(n,k,d,m)\right| = \left| B(n,k,d,m)\right|$.
\end{theorem}
\noindent Of course, for $m\rightarrow\infty$ (indeed, for $m > n$), this result reduces to Theorem \ref{sythm}.

We give two different proofs of Theorem \ref{stheorem}.  First, we demonstrate the following identity, a $q$-series extension of (\ref{solutionI}) used in Madhyastha's proof:
\begin{align}
\frac{q^{dk}}{(q^d;q^d)_k} \frac{(q^d;q^d)_{m}}{(q;q)_{dm}} =\begin{cases}
\frac{q^{kd}}{(q^{d(m+1)};q^d)_{k-m}}\frac{1}{(q;q)_{md}}, & \text{ if } m < k,\\
& \\
\frac{q^{kd}}{(q;q)_k}\frac{(q^{d(k+1)};q^d)_{m-k}}{(q^{k+1};q)_{m-k}}\frac{1}{(q^{m+1};q)_{md-m}}, & \text{ if } m\ge k.
\end{cases}\label{finF}
\end{align}
We first interpret (\ref{finF}).  Examining the left-hand side,
\begin{align}
\frac{q^{dk}}{(q^d;q^d)_k} \frac{(q^d;q^d)_{m}}{(q;q)_{dm}},\label{finFA1}
\end{align}  we notice that the first factor of (\ref{finFA1}) enumerates partitions into exactly $k$ parts divisible by $d$.  The second factor of (\ref{finFA1}) enumerates partitions in which \textit{no} parts are divisible by $d$, but moreover one in which all parts are bounded above by $md$; indeed, because the part $md$ is specifically excluded, the parts must be \textit{strictly less} than $md$.  Thus the left-hand side of (\ref{finF}) enumerates partitions of $A(n,k,d,m)$.

Examining the right-hand side of (\ref{finF}) in the case that $m<k$, we note that the very last factor, $1/(q;q)_{md}$, enumerates partitions into parts of size up to $md$, while $1/(q^{d(m+1)};q^d)_{k-m}$ enumerates partitions of size larger than $md$ but smaller than $kd$, in which the parts are divisible by $d$.  Finally, the factor of $q^{kd}$ confirms that $kd$ does in fact occur.  Thus, we have partitions of type $B(n,k,d,m)$ for $m<k$ enumerated by the right-hand side of (\ref{finF}).

For $m\ge k$, $1/(q;q)_k$ enumerates partitions into parts of size up to and including $k$ with no additional restrictions.  To insure that $k$ occurs at least $d$ times, we have the factor $q^{kd} = q^{k+k+\cdots+k}.$  On the other hand, the factor $1/(q^{m+1};q)_{md-m}$ enumerates partitions into parts of size strictly greater than $k$ and less than or equal to $md$, but with no other restrictions.
Finally, the factor
\begin{align*}
\frac{(q^{d(k+1)};q^d)_{m-k}}{(q^{k+1};q)_{m-k}}
\end{align*} consists of factors of the form
\begin{align*}
\frac{1-q^{dj}}{1-q^j} = 1+q^j+q^{2j}+\cdots+q^{(d-1)j},
\end{align*} for $k < j\le m$.  Thus, the parts strictly larger than $k$ and less than or equal to $m$ must occur fewer than $d$ times.  Thus the right-hand side of (\ref{finF}) enumerates $B(n,k,d,m)$ in both cases.

\begin{proof}
We have two cases to prove; but in both cases, we begin with (\ref{finFA1}).  First we suppose that $m<k$.  Then we have
\begin{align*}
(q^d;q^d)_k = (q^d;q^d)_m (q^{d(m+1)};q^d)_{k-m}.
\end{align*}  Thus substituting this for $(q^d;q^d)_k$ in (\ref{finFA1}), we have
\begin{align*}
\frac{q^{dk}}{(q^d;q^d)_k} \frac{(q^d;q^d)_{m}}{(q;q)_{dm}} &= 
\frac{q^{dk}}{(q^d;q^d)_m (q^{d(m+1)};q^d)_{k-m}} \frac{(q^d;q^d)_{m}}{(q;q)_{dm}}\\
&= \frac{q^{dk}}{(q^{d(m+1)};q^d)_{k-m}} \frac{1}{(q;q)_{dm}}.
\end{align*}
Similarly, for $m\ge k$,
\begin{align*}
(q^d;q^d)_m = (q^d;q^d)_k (q^{d(k+1)};q^d)_{m-k}.
\end{align*}  Substituting into (\ref{finFA1}), 
\begin{align}
\frac{q^{dk}}{(q^d;q^d)_k} \frac{(q^d;q^d)_{m}}{(q;q)_{dm}} &= 
\frac{q^{dk}}{(q^d;q^d)_k} \frac{(q^d;q^d)_k (q^{d(k+1)};q^d)_{m-k}}{(q;q)_{dm}}\label{finRm1}\\
&= q^{dk}\frac{(q^{d(k+1)};q^d)_{m-k}}{(q;q)_{dm}}.\label{finRm2}
\end{align}  Next, we note that
\begin{align*}
(q;q)_{dm} = (q;q)_k (q^{(k+1)};q)_{m-k}(q^{(m+1)};q)_{md-m}.
\end{align*}  With this substitution into (\ref{finRm1})--(\ref{finRm2}) we have 
\begin{align*}
q^{dk}\frac{(q^{d(k+1)};q^d)_{m-k}}{(q;q)_{dm}} &= q^{dk}\frac{(q^{d(k+1)};q^d)_{m-k}}{(q;q)_k (q^{(k+1)};q)_{m-k}(q^{(m+1)};q)_{md-m}}\\
&= \frac{q^{dk}}{(q;q)_k}\frac{(q^{d(k+1)};q^d)_{m-k}}{(q^{(k+1)};q)_{m-k}}\frac{1}{(q^{(m+1)};q)_{md-m}}.\hfill\qedhere
\end{align*}\end{proof}
However, we have also discovered a corresponding extension of the direct constructive proof that we already gave.  We note that our new proof depends on the proof of a finite version of Glaisher's theorem which can be found in \cite{Nyirenda}; the latter is itself an extension of the finite version of Euler's partition identity proved by Andrews in \cite{Andrewsf}.  Both of these finite bound theorems play a substantial role in achieving finite-bound equinumerosity results, e.g., \cite{Fu}.

\begin{proof}

Let $\lambda\in A(n,k,d,m)$.  We write
\begin{align*}
\lambda = \mu + o,
\end{align*} in which $\mu$ contains the $k$ parts divisible by $d$, and $o$ contains the remaining parts.

We begin by considering the subpartition $\mu$.  Here we map $\mu$ to its conjugate, $\mu^{\ast}$, in which $k$ occurs at least $d$ times (indeed, the number of occurrences of $k$ is a positive multiple of $d$).  We can therefore write $\mu^{\ast}$ as
\begin{align*}
\mu^{\ast}: \sum_{i=1}^k dq_i\cdot (i),
\end{align*}  in which we interpret the parts as $i$, and the number of occurrences as $dq_i$ for some $q_i\in\mathbb{Z}_{\ge 0}$ (with $q_k > 0$).

We will in turn break $\mu^{\ast} = \mu^{\ast}_0+\epsilon$ into two smaller subpartitions:
\begin{align*}
\mu^{\ast}_0 &: \sum_{i=1}^{\mathrm{min}(m,k)} dq_i \cdot (i),\ \ \ \ \ \epsilon: \sum_{m< i\le k} q_i \cdot (di).
\end{align*}  For $\epsilon$ we reinterpret the terms: the parts are now $di$ as $i$ ranges from 1 to $m$.  Notice that for $m < k$, the largest part from either of these subpartitions will be $dk$, and it will come from $\epsilon$.  Note also that $\epsilon$ is empty for $m\ge k$.

Next, we consider the subpartition $o$.  Notice that these parts are bounded strictly above by $md$, and that none of them are divisible by $d$.  Write
\begin{align*}
o: \sum_{\substack{1\le j< md,\\ d\nmid j}} N_j\cdot (j),
\end{align*} in which $j$ is a given part, and $N_j\ge 0$ the number of occurrences of that part.

For each $j\le md$ there exists a unique integer $L_j\ge 0$ such that
\begin{align*}
m< jd^{L_j} \le md.
\end{align*}  We write $N_j$ in base $d$, but we collect together all terms corresponding to powers of $d$ which match or exceed $L_j$:
\begin{align}
N_j &= \sum_{l\ge 0} a_{j,l}\cdot d^l = \sum_{0\le l\le L_j-1} a_{j,l}\cdot d^l + M_j\cdot d^{L_j}.\label{omapdeltaA}
\end{align}  Here for $l\ge 0$, $0\le l\le L_j-1$ we have $0\le a_{j,l}\le d-1$.  Now we re-examine our original subpartition $o$:
\begin{align*}
&\sum_{\substack{1\le j< md,\\ d\nmid j}} N_j\cdot (j) =\sum_{\substack{1\le j< md,\\ d\nmid j}}\left( \sum_{0\le l\le L_j-1} a_{j,l}\cdot d^l + M_j\cdot d^{L_j} \right)\cdot (j).
\end{align*}  We thus have our new subpartition
\begin{align}
\delta: \sum_{\substack{1\le j< md,\\ d\nmid j,\\ 0\le l\le L_j-1}} a_{j,l}\cdot \left(j\cdot d^l \right) + \sum_{\substack{1\le j< md,\\ d\nmid j}} M_j\cdot \left(j\cdot d^{L_j}\right).\label{omapdeltaZ}
\end{align}  The first sum contributes parts of the form $i=j\cdot d^l$ for $0\le l\le L_j - 1$, and therefore $1\le i\le m$; indeed, by unique factorization, any possible part between 1 and $m$ must have this form.  Moreover, each $i$ occurs less than $d$ times.  

The parts which are permitted to repeat arbitrarily many times are $i=j\cdot d^{L_j}$, with $m < i\le md$.  Again, by the possible range of $j$ and through unique factorization, every possible part between $m$ and $md$ must have this form.  These parts repeat $M_j$ times, for any integer $M_j\ge 0$.

We will send $\lambda$ to the following partition:
\begin{align*}
\kappa = \mu^{\ast}_0 + \epsilon + \delta.
\end{align*}  We need to show that $\kappa\in B(n,k,d,m)$.\\

\noindent\textbf{Case: $m < k$}\\
In this case, $\epsilon$ is not empty, and has largest part $dk$.  The parts of $\mu^{\ast}_0$ are all $\le m < k$, and the parts of $\delta$ are all $\le md < kd$.  Thus our largest part will be $kd$.

A part $i$ with $1\le i\le m$ will have the form $i=jd^l$ for some $j$ between 1 and $m$, and some nonnegative power $l$ of $d$.  It occurs $dq_{i} + a_{j,l}$ times; here $dq_{i}$ comes from $\mu^{\ast}_0$, and $a_{j,l}$ comes from $\delta$.  Per the conditions of the division algorithm, this is an unrestricted nonnegative integer.

Any part $i$ with $m< i\le md$, has the form $jd^{L_j}$ and also occurs an unrestricted number $M_j$ of times, coming from $\delta$.  Neither $\mu^{\ast}_0$ nor $\epsilon$ contribute to these parts.

Finally, for parts greater than $md$, we have only terms $dj$ from $\epsilon$.  The only restrictions is by size (the parts must be greater than $md$), and divisibility by $d$.

We thus have a well-defined partition $\kappa\in B(n,k,d,m)$ with no set restrictions beyond those of $B(n,k,d,m)$.\\

\noindent\textbf{Case: $m\ge k$}\\
In this case, $\epsilon$ is empty, and the largest part from $\mu^{\ast}_0$ is $k$.  Indeed, any part $i$ with $1\le i\le k$ occurs $dq_i$ times from $\mu^{\ast}_0$.  It may also occur in $\delta$; however, because $i\le k \le m$, it has the form $i=jd^l$ with $d\nmid j$ and $0\le l\le L_j-1$, and can only occur $a_{j,l}$ times, i.e., no more than $d-1$ times.  Thus the number of occurrences of $i$ in $\kappa$ has the form $dq_i+a_{j,l}$; it is therefore arbitrary, per the conditions of the division algorithm.  Specifically for $i=k$, we have the additional condition that $q_k > 0$, so that $k$ must occur at least $d$ times.

Similary, all parts $i > m$ can only come from $\delta$, and be $\le md$.  Moreover, we have already shown that such an $i$ may repeat $M_j$ times, with $i=jd^{L_j}$, $d\nmid j$.

Only the parts $i$ with $k < i \le m$ have a restriction on their occurrences, appearing less than $d$ times.

Finally, the largest parts in $\mu^{\ast}$ are $k\le m\le md$, and the largest possible parts of $\delta$ are $md$.  Thus $\kappa\in B(n,k,d,m)$.\hfill\qedhere\end{proof}

For the sake of illustration, we give two examples.  First, we take $d=3$, $k=7$, $m=4$, $n=123$.  Notice that $m<k$.  Moreover, we have $md=12$ and $kd=21$.  We consider the partition
\begin{align*}
\lambda &:= (2)15+12+11+9+8+(4)7+(2)6+5+3+(2)2+1.
\end{align*}  The parts divisible by 3 are 15 (occuring twice), 12, 9, 6 (occuring twice), and 3.  That is, there are $k=7$ parts divisible by $d=3$.  The largest part not divisible by 3 is $11<12=md$.  Thus, $\lambda\in A(123,7,3,4)$.

We build our bijection first by separating $\lambda$ into its key subpartitions:
\begin{align*}
\mu:\ & (2)15+12+9+(2)6+3,\\
o:\ & 11+8+(4)7+5+(2)2+1.
\end{align*}  First,we map $\mu$ to its conjugate $\mu^{\ast}$:
\begin{align*}
& (2)15+12+9+(2)6+3\longleftrightarrow (3)7+(3)6+(3)4+(3)3+(3)2.
\end{align*}  Notice that each part $i$ of $\mu^{\ast}$ lies between 1 and $k=7$, and occurs a multiple of 3 times (possibly 0).

In turn we split $\mu^{\ast}$ into $\mu^{\ast}_0$ and $\epsilon$.  Because $\mathrm{min}(m,k)=m$, $\epsilon$ is not empty.  We construct $\epsilon$ by examining the parts of $\mu^{\ast}$ which are greater than $m=4$.  The factor of 3 in the number of occurences is now multiplied by our part to produce a new part.  What is left of the number of occurences after dividing out 3 is our new occurence:
\begin{align*}
(3)7+(3)6 &= (3\cdot 1)7+(3\cdot 1)6 \longleftrightarrow (1)21+(1)18 =\epsilon.
\end{align*}  The remaining parts of $\mu^{\ast}$ give us $\mu^{\ast}_0$:
\begin{align*}
\mu^{\ast}_0 = (3)4+(3)3+(3)2.
\end{align*}
We now turn our attention to
\begin{align*}
o:\ & 11+8+(4)7+5+(2)2+1.
\end{align*}  For $1\le j< 12, 3\nmid j$, there exists a unique $L_j$ such that $4 <jd^{L_j}\le 12$.  We may compute $L_j$ for each part $j$ of $o$:
\begin{align*}
&\text{For } j=11,8,7,5, L_j=0.\\
&\text{For } j=2, L_j=1.\\
&\text{For } j=1, L_j=2.
\end{align*} Thus for $j=1,2$, we can expand the number of occurences by base $3$ up to the power of 2 for $j=1$, and $1$ for $j=2$.  But the occurences for these parts in $o$ are 1 and 2, respectively, with resulting maximum power being 0 in the base $3$ expansions.  Therefore we can simply apply Glaisher's bijection.  In this case, it is a trivial mapping:
\begin{align*}
(2)2+(1)1\longrightarrow (2)2+(1)1.
\end{align*} Since $L_j=0$ for the remaining parts of $o$, we will not apply Glaisher's bijection at all; we should multiply each $j$ by $3^{L_j}$.  However, this factor is simply 1, and we may therefore simply reproduce the partition $o$:
\begin{align*}
\delta:\ & 11+8+(4)7+5+(2)2+1.
\end{align*}  Thus we have
\begin{align*}
\kappa &= (1)21+(1)18 + (3)4+(3)3+(3)2 + 11+8+(4)7+5+(2)2+1\\
&= 21+18+11+8+(4)7+5+(3)4+(3)3+(5)2+1.
\end{align*}  Our largest part is $kd=(7)(3)=21$, and the parts greater than $md=12$ are 21, 18, both divisible by $d=3$.  Thus $\kappa\in B(123,7,3,4)$.

For our second example, we consider the $d=3$, $k=4$, $m=7$, $n=189$.  Notice that $m\ge k$.  Moreover, we have $md=21$ and $kd=12$.  Let us take, say,
\begin{align*}
\lambda := 24+21+20+17+15+(4)14+9+(2)7+(5)2+(3)1.
\end{align*}  Here, the parts divisible by 3 are 24, 21, 15, 9.  Thus there are $k=4$ parts divisible by $d=3$.  The largest part not divisibly by 3 is 20, which is less than $md=21$.  Thus, $\lambda\in A(189,4,3,7)$.  Again, we split $\lambda$ into 
\begin{align*}
\mu &:\ 24+21+15+9,\\
o &:\ 20+17+(4)14+(2)7+(5)2+(3)1.
\end{align*}  We map $\mu$ to its conjugate $\mu^{\ast}$:
\begin{align*}
24+21+15+9\longleftrightarrow (9)4+(6)3+(6)2+(3)1.
\end{align*}  Because $m>k$, the subpartition $\epsilon$ is empty, and $\mu^{\ast}_0 = \mu^{\ast}$.

Next we take
\begin{align*}
o:\ & 20+17+(4)14+(2)7+(5)2+(3)1.
\end{align*} For $1\le j< md, d\nmid j$, we compute $L_j$ as before, and we find that
\begin{align*}
&\text{For } j=20,17,14, L_j=0.\\
&\text{For } j=7, L_j=1.\\
&\text{For } j=2,1, L_j=2.
\end{align*}  Now the mapping is nontrivial, and using (\ref{omapdeltaA})--(\ref{omapdeltaZ}) we get
\begin{align*}
&[(1)20+(1)17+(4)14]+[(2)7+(2\cdot 3^0+1\cdot 3^1)2+(1\cdot 3^1)1]\\
=&[(1)20+(1)17+(4)14]+[(2)7+(2)2+(1)6+(1)3].
\end{align*}
Putting our new subpartitions back together, we get
\begin{align*}
&[(9)4+(6)3+(6)2+(3)1]+[(2)7+(2)2+(1)6+(1)3]\\ &+[(1)20+(1)17+(4)14],\\
\kappa &= 20+17+(4)14+(2)7+6+(9)4+(7)3+(8)2+(3)1.
\end{align*} Notice that $k=4$ occurs at least $d=3$ times (indeed, it occurs 9 times), and all parts are bounded above by $md=21$.  Of particular interest are the parts larger than $k=4$ and less than or equal to $m=7$: $(2)7+6$.  These parts occur less than 3 times.  Thus, $\kappa\in B(189,4,3,7)$.

\end{document}